\documentclass[amsthm]{elsart}

\usepackage{yjsco}
\usepackage{natbib}
\usepackage{amsmath,amssymb}
\usepackage{amssymb}
\usepackage{amsfonts}
\usepackage{graphicx}
\usepackage{epsfig}
\usepackage{color}
\usepackage{url}

\newcommand{\eq}{\Leftrightarrow}

\newcommand{\NN}{\mathbb{N}}
\newcommand{\ZZ}{\mathbb{Z}}

\newtheorem{theorem}{Theorem}[section]

\newtheorem{proposition}[theorem]{Proposition}
\newtheorem{corollary}[theorem]{Corollary}

\begin{document}

\begin{frontmatter}

\title{$G$-graphs Characterisation and Incidence Graphs}

%\thanks{This research was partly supported by .....}

\author{David Ellison}
\address{RMIT University, 124 Little La Trobe St, Melbourne VIC 3000, Australia}
\ead{e-mail: davidellison@polytechnique.edu}
%\ead[url]{URL 1}

\author{Ruxandra Marinescu-Ghemeci}
\address{University of Bucharest, Str. Academiei nr.14, Bucharest, Romania}
\ead{e-mail: verman@fmi.unibuc.ro}
%\ead[url]{URL 2}

\author{Cerasela Tanasescu}
\address{ESSEC Business School, Av. B. Hirsch, 95000 Cergy Pontoise, France\\CEREGMIA, University of Antilles-Guyane, Schoelcher, Martinique, France}
\ead{e-mail:  tanasescu@essec.edu}
%\ead[url]{URL 2}

\begin{abstract}
Graphs derived from groups are a widely studied class of graphs, motivated by their highly symmetric structure. In particular, $G$-graphs offer an easy and interesting alternative construction of semi-symmetric graphs. After recalling the main properties of these graphs, this papers gives an extended characterisation of $G$-graphs and develops the link between bipartite $G$-graphs and incidence graphs. It appears that these two classes of graphs have a wide overlapping despite having completely different constructions. We give partial answers to the problem of finding which complete simple graphs have a $G$-graph as their incidence graph.

\end{abstract}


\end{frontmatter}
\section{Introduction}\label{sec:intro}
This paper is a contribution to the study of graphs derived from algebraic groups. The most popular graphs defined by a group are Cayley graphs. $G$-graphs correspond to an alternative construction. These graphs, introduced in \cite{brfa0}, have highly regular properties like Cayley graphs. In particular, because the algorithm for constructing $G$-graphs is simple, it is a useful tool for constructing new symmetric and semi-symmetric graphs (see \cite{bre3,bre2}). One interesting direction is  the study of the properties of $G$-graphs and the characterisation theorems for $G$-graphs. In the following sections, we propose such a characterisation and study the incidence graphs of $G$-graphs.
The paper is organised as follows. In Section 2 we first recall the main definitions from group and graph theories and we present the class of $G$-graphs with some of their basic properties. In Section 3, we propose a generalisation of the characterisation theorem. Based on these results, we investigate in Section 4 the connected components of $G$-graphs as well as the properties of complete bipartite $G$-graphs. More precisely, we give the proof that complete bipartite multigraphs without loops are abelian $G$-graphs. Section 5 is entirely dedicated to the study of incidence graphs. First, we define the incidence graphs and recall their main properties. In the second part we establish a link between incidence graphs and bipartite  $G$-graphs. We also give a necessary condition and a sufficient condition for the incidence graph of a bipartite $G$-graph to be a $G$-graph. In the third part of this section, we study the case of the incidence graphs of complete graphs.

\section{Basic Definitions}~\label{sec: def}
This section contains basic definitions of algebraic graph theory. Further definitions may be found in \cite{rob} and \cite{ref_west}.

\subsection{Groups}
\noindent Let $(G,\cdot,e)$ be a  group, $e$ denotes the identity element of $G$ and "$\cdot$" denotes a canonical operation (multiplicative notations). For every $g$ in $G$, the set $\{g^n|n \in \mathbb{Z}\}$  forms a subgroup of $G$ called the {\em cyclic group generated by $g$}. We denote it by $\langle g \rangle$. Its cardinal is called the {\em order} of $g$.

\bigskip\noindent Let $S$ be a set of elements in $G$. The {\em subgroup generated by $S$}, denoted by $\langle S\rangle$, is the smallest subgroup of $G$ which contains $S$. If $S=\{s_1,...,s_n\}$, $\langle S\rangle$ can also be denoted by $\langle s_1,...,s_n\rangle$. If $\langle S\rangle=G$, $S$ is said to be a {\em generating set} of $G$.

\bigskip\noindent Given two groups $G$ and $H$, a {\em homomorphism} from $G$ to $H$ is a function $\varphi:G\longrightarrow H$ such that $\varphi(gg')=\varphi(g)\varphi(g')$, for every $g,g'\in G$. A bijective homomorphism from $G$ to $H$ is called {\em isomorphism}. An isomorphism from $G$ onto itself is called {\em automorphism}. We denote by $Aut(G)$ the group of automorphisms of $G$ under composition law.

\bigskip\noindent We define a {\em(left) group action} of $G$ on  a set $X$ as a function from $G\times X$ to $X$, ${\left(g, x\right) \rightarrow g\cdot x}$ satisfying conditions: $e\cdot x = x$ for every $x\in X$ and $g\cdot \left(g^{'}\cdot x\right) = \left(g\cdot g^{'}\right)\cdot x$ for every $g$,  $g^{'}\in G$ and $x\in X$. The {\em orbit} of $x$, noted $O_x$, is the set $\{g\cdot x,g\in G\}$. The orbits of the action of $G$ on $X$ form a partition of $X$. The action is {\em transitive} if $\forall x,y \in X, \exists g \in G$ such that $g\cdot x = y$.  For $x\in X$ we define the {\em stabilizer subgroup} of $x$ as the set of all elements in $G$ that fix $x$: $Stab_{G}x = \left\{g : g \cdot x = x\right\}$. The action is said to be {\em regular} if for all $(x,y)\in X^2$, there is a unique $g\in G$ such that: $g\cdot x=y$.

\bigskip\noindent For instance, if $H$ is a fixed subgroup of a group $G$, $H$ acts on $G$ via the group law. For $x\in G$, the orbit of $x$ is $Hx=\{hx |h \in H\}\subset G$ and is called {\em right coset of H containing x}. It follows that the right cosets of $H$ form a partition of $G$: for any $x,y\in G$, $Hx=Hy$ or $Hx\cap Hy=\emptyset$.

\bigskip\noindent The axiom of choice is an axiom of set theory equivalent to the statement that the cartesian product of a collection of non-empty sets is non-empty. More precisely it states that for every indexed family $(S_i)_{i \in I}$ of nonempty sets there exists an indexed family $(x_i)_{i \in I}$ of elements such that $x_i \in S_i$ for every $ i \in I$.  A consequence of the axiom of choice which will be of use is: every set is in one-to-one correspondence with an abelian group (see \cite{chapeau}).

\subsection{Graphs}

%Bondy
\noindent An non directed multigraph $\Gamma$ is a triple $(V(\Gamma),E(\Gamma), \psi_\Gamma)$, where $V(\Gamma)$ is the set of vertices,  $E(\Gamma)$which is disjoint from $V(\Gamma)$ is the set of edges and $\psi_\Gamma$ is an incidence function  that associates to each edge of $\Gamma$ an unordered pair of (not necessarily distinct) vertices of $\Gamma$.  If $\psi_{\Gamma}(e)=\{u,v\}$, then the edge $e$ is said to {\em link} vertices $u$ and $v$. For any two adjacent vertices $u$ and $v$, the {\em multi-edge} $uv$ is the set of edges with end-points $u$ and $v$.

\bigskip\noindent For notation simplicity, if no confusion occurs, we will write $V$, $E$ and $\psi$ instead of $V(\Gamma)$, $E(\Gamma)$ and $\psi_\Gamma$.

\bigskip\noindent An edge with identical end-points is called a {\em loop}. Two or more edges with the same pair of end-points are said to be {\em parallel edges}. A multigraph is {\em simple} if it has no loops or parallel edges. In any simple graph, we may dispense with the incidence function $\psi$ by renaming each edge as the unordered pair of its end-points. 

%\begin{rem}\label{rem:multigraph}
%Sometimes multigraphs are simply called {\em graphs}. However, this may cause some confusion with the most usual terminology. That is why we will specify multigraph when explicitly working with parallel edges.
%\end{rem}

\bigskip\noindent The {\em degree} of a vertex $v$ in a multigraph $\Gamma$, denoted by $d_{\Gamma}(v)$ or simply $d(x)$ when no ambiguity occurs is the number of edges of $\Gamma$ incident with $v$, each loop counting as two edges. 

\bigskip\noindent A {\em finite walk} is an alternating sequence of vertices and edges $(v_0,a_1,v_1,...,a_n,v_n)$ such that for all $i$, $v_{i-1}$ and $v_i$ are the end-points of $a_i$. $n$ is called the {\em length} of the walk. A {\em cycle} or {\em n-cycle} is a walk $(v_0,a_1,v_1,...,a_n,v_n)$ where $v_0=v_n$ and no other two vertices in the walk are the same.

\bigskip\noindent Let $\Gamma=(V,E,\psi)$ and let $X\subset V$. The {\em induced subgraph} $\Gamma[X]$ is the graph $(X,E',\psi|_{E'})$ with $E'=\psi^{-1}(P_2(X))$, where $P_2(X)$ is the set of pairs of (not necessarily distinct) elements of $X$.

\bigskip\noindent A {\em clique} of $\Gamma$ is a set of vertices of $\Gamma$ which are all adjacent. A {\em stable} is a set of vertices no two of which are adjacent.

\bigskip\noindent Let $\Gamma_{1}$ and $\Gamma_{2}$ be two multigraphs. A {\em graph homomorphism} is a couple $(f,f^\#)$ where $f: V(\Gamma_{1})\longrightarrow  V(\Gamma_{2})$ and $f^{\#} : E(\Gamma_1) \longrightarrow E(\Gamma_2)$ are such that $\psi_{\Gamma_1}(a) = \{u,v\}$ if and only if  $\psi_{\Gamma_2}(f^{\#}(a)) = \{f(u),f(v)\}$. A {\em graph isomorphism} is a graph homomorphism $(f,f^\#)$ where $f$ and $f^\#$ are bijective. $\Gamma_{1}$ and $\Gamma_{2}$ are {\em isomorphic}, written $\Gamma_1\cong\Gamma_2$, if there is a graph isomorphism between the two. An isomorphism from a graph $\Gamma$ onto itself is called an {\em automorphism} of $\Gamma$. $Aut(\Gamma)$ denotes the group of automorphisms of a graph $\Gamma$ under composition law. The identity of $Aut(\Gamma)$ is written $Id_\Gamma$.

\bigskip\noindent In the case of simple graphs, the definition of isomorphism can be stated more concisely.  If  $(f, f^{\#})$ is an isomorphism between simple graphs $\Gamma_{1}$ and $\Gamma_{2}$, then  $f^{\#}$ is completely determined by $f$. Consequently, an isomorphism between two simple graphs $\Gamma_{1}$ and $\Gamma_{2}$ as a bijection $f: V(\Gamma_{1})\longrightarrow  V(\Gamma_{2})$  which preserves adjacency.

% The actions of the group $Aut(\Gamma)$ on the set $V_\Gamma$ and $E_\Gamma$ are defined as follows: $(f,f^\#)\cdot v=f(v)$ and $(f,f^\#)\cdot a=f^\#(a)$, for all $f\in Aut(\Gamma), v\in V_\Gamma,$ and $a\in E_\Gamma$.

\subsection{$G$-Graphs}

\noindent The  literature on $G$-graphs (sometimes called orbit graphs) already includes several definitions, which produce slightly different objects. The first definition of $G$-graphs was given in \cite{brfa0}. In the original definition, the vertices of $G$-graphs were cycles of the left $s$-translation on $G$ $(x,sx,...,s^{o(s)-1}x)$. The $G$-graph was defined as a multigraph with a $p$-edge between any to cycles with an intersection of $p$ elements. Also, there were $o(s)$ loops on the vertex $(x,sx,...,s^{o(s)-1}x)$. Using the same definition, \cite{breCayley}, \cite{bre3}, \cite{symmetry} and \cite{bre2} also introduce $G$-graphs without loops. \cite{orbitgraphs} introduces orbit graphs, which are slightly different from the original $G$-graphs: the first difference is that the vertices are seen as orbits of an action instead of cycles. Also, the definition of orbit graphs includes a colouring (labeling) of the edges. In \cite{balcor}, the vertices are defined as right cosets of cyclic groups. The main difference is that $G$-graphs are intersection graphs instead of multigraphs.

\bigskip\noindent In this paper, we chose to define $G$-graphs as multigraphs without loops and with a labeling. We introduce a novelty by using multisets, allowing repetitions of the generating elements. The reason for these choices is that we obtain simpler and more general versions of the main properties, in particular of the $G$-graphs characterisation.

\bigskip\noindent A {\em multiset} $X$ is a usual generalisation of sets where each element may have several occurrences. 

\bigskip\noindent Let $G$ be a group and $S$ a multiset of elements of $G$. The {\em $G$-graph} $\Phi(G,S)$ is the labeled multigraph $(V,E,\psi)$ such that:
\begin{enumerate}
\item The set of vertices of $\Phi (G,S)$ is $V = \bigcup_{s \in S}V_{s}$, where $V_{s}=\{\langle s \rangle x, x\in G\}$ ($\langle s \rangle x$ is the right coset of $\langle s \rangle$ containing $x$). . 
\item For  $\langle s\rangle x, \langle t \rangle y\in V$ ($s\neq t$), there exists an edge between $\langle s \rangle x$ and $\langle t \rangle y$ labeled g for each $g \in \langle s \rangle x\; \cap \langle t \rangle y$. 
\end{enumerate}

\bigskip\noindent The edge between $\langle s \rangle x$ and $\langle t \rangle y$ labeled $g$ will be noted $(\langle s\rangle x\langle t \rangle y,g)$.

\bigskip\noindent The {\em $G$-graph with loops} $\Psi(G,S)$ is the graph obtained by adding the corresponding loops to $\Phi(G,S)$. 

\bigskip\noindent A repetition of an element $s$ in the multiset $S$ induces a repetition of the right cosets of $\langle s \rangle$ in $V$. For purely formal purposes, in order for $V$ to be a set rather than a multiset, it may be necessary to distinguish between such repetitions by denoting the vertices by $(\langle s \rangle x,k)$, where $k$ is an index belonging to the set of occurrences of $s$ in $S$. However, instead of using this overburdened notation, we will simply allow repetitions in $V$.

\bigskip\noindent Given a $G$-graph $\Gamma=\Phi(G,S)$, and $s\in S$, the {\em level of $s$}, noted $V_s$, is the stable of $\Gamma$ which comprises all the vertices of the form $\langle s \rangle  x, x\in G$. The elements of a level of $\Gamma$ form a partition of $G$: it is the partition of $G$ into right cosets of $\langle s \rangle $.

\bigskip\noindent Given a $G$-graph $\Gamma=\Phi(G,S)$, and $g\in G$, the {\em colour clique} of $g$, denoted by $C_g$, is the set of vertices of $\Gamma$ which contain $g$.

\bigskip\noindent Given a $G$-graph $\Gamma=\Phi(G,S)$, and $g\in G$, the {\em shift} $\delta_g$ is the automorphism of $\Gamma$ defined by $\delta_g(\langle s \rangle x)=\langle s \rangle  xg$ and $\delta_g^\#(\langle s \rangle x\langle t \rangle  y,h)=(\langle s \rangle xg\langle t\rangle yg,hg)$.

\bigskip\noindent The following proposition presents 5 key elements of the structure of $G$-graphs. The first two points can be found in \cite{breCayley}, \cite{bre3}, \cite{orbitgraphs}, \cite{symmetry} and \cite{bre2}. We generalise previous results to the infinite case and to the case where there are repetitions in $S$.

\begin{proposition}\label{structure}
Let $\Phi(G,S)$ be a $G$-graph, 
\begin{enumerate} 
\item $G\simeq \{\delta_g,g\in G\}$
\item $\forall g\in G,\forall s\in S, \delta_g(V_s)=V_s$
\item $\forall (g,g')\in G^2, \delta_{g'}(C_g)=C_{gg'}$
\item $\forall (s,s')\in S^2, s\neq s'$, $E({\Psi(G,S)[V_s]})$ and $E({\Phi(G,S)[V_s\cup V_{s'}]})$ are in one-to-one correspondence with $G$, and $E({\Psi(G,S)})$ is in one-to-one correspondence with $G\times P_2(S)$.
\item For any pair of adjacent vertices $\langle s \rangle x$ and $\langle t \rangle y$, the set of labels of the edges of the multi-edge $\langle s \rangle x \langle t \rangle y$ is a right coset of $\langle s \rangle \cap\langle t \rangle $.
\end{enumerate}
\end{proposition}

\noindent\textbf{Proof:}
(1) $\{\delta_g,g\in G\}$ is a group under the composition law. $g\mapsto\delta_g$ is a group homomorphism, and its inverse is obtained by associating with $\delta_g$ the label of the image of any edge labeled $e$ in $\Phi(G,S)$.

\bigskip\noindent (2) Let $s\in S$ and $(g,x)\in G^2$. $\delta_g(\langle s \rangle x)=\langle s \rangle xg\in V_s$ and $\langle s \rangle x=\delta_g(\langle s \rangle xg^{-1})\in\delta_g(V_s)$. Hence $\delta_g(V_s)=V_s$.

\bigskip\noindent (3) $C_g=\{\langle s \rangle g,s\in S\}$. Hence $\delta_{g'}(C_g)=\{\langle s \rangle gg',s\in S\}=C_{gg'}$.

\bigskip\noindent (4) Let $s,t\in S$ and $g\in G$. Since the elements of each level form a partition of $G$, there is exactly one vertex in both $V_s$ and $V_{t}$ which contains $g$ (these vertices are $\langle s\rangle g$ and $\langle s'\rangle g$). Hence, there is exactly one element of $V_s$ with a loop labeled $g$, and exactly one link labeled $g$ joining $V_s$ and $V_{t}$. Therefore, the map which associates with each edge the levels of its end-points and its label defines a one-to-one correspondence between $E({\Psi(G,S)})$ and $G\times P_2(S)$.

\bigskip\noindent (5) Let $h\in \langle s \rangle x\cap\langle t \rangle y$.\\
$z\in\langle s \rangle x\cap\langle t \rangle y\eq \exists (m,n)\in\ZZ^2:z=s^mh$ and $z=t^nh$\\
$\eq \exists w\in\langle s \rangle \cap\langle t \rangle:z=wh$\\
Hence $\langle s \rangle x\cap\langle t \rangle y=(\langle s \rangle \cap\langle t \rangle)h$

\bigskip\noindent This result shows us that $G$-graphs can be decomposed following two transversal directions. The set of vertices can be partitioned either into levels, which are stables, or into colour cliques. The vertex $\langle s \rangle g$ is the intersetion of $V_s$ and $C_g$. Each shift sends one colour clique to another while stabilising the levels. A direct consequence of (5) is that $\Phi(G,S)$ is simple if and only if the elements of $S$ are pairwise independent.

\section{Characterisations of $G$-graphs with loops and $G$-graphs}

A first partial characterisation of $G$-graphs in the bipartite case is given in \cite{bre3}, Th 6.1, and \cite{symmetry}, Th 4.1. Two slightly different extensions of this theorem can be found in \cite{orbitgraphs}, Th 1, and in \cite{balcor}, Th 2. We further extend the characterisation of $G$-graphs in two different ways: by allowing repetitions in $S$ and by including the infinite case. We are thus able to recognise a wider class of graphs, now including simple complete graphs.

\bigskip In this article, we present two versions of the characterisation: for $G$-graphs with loops and without loops. Adding loops makes the characterisation simpler, as it requires one fewer condition. The proof also becomes more efficient. Yet, we will use the version without loops in the applications.

\begin{theorem}\label{characterisation} Characterisation of $G$-graphs with loops: \\
A multigraph $\Gamma=(V,E,\psi)$ is a $G$-graph with loops if and only if there is a subgroup $H$ of $Aut(\Gamma)$ and a clique $C$ with loops such that:
\begin{enumerate}
\item $C$ intersects every orbit of the action of $H$ on $V$
\item $\forall u\in C, Stab_Hu$ is cyclic
\item $\forall u\in C$, for each orbit $O$ of $H$, $Stab_Hu$ acts regularly on the set of edges adjacent to $u$ with their other end-point in $O$.
\end{enumerate}
\end{theorem}

\noindent\textbf{Proof:}
Let $\Gamma=\Psi(G,S)=(V,E,\psi)$, we will show that $H=\{\delta_g,g\in G\}$ and $C_e$ verify the above.

\bigskip\noindent (1) It stems from proposition~\ref{structure} that the orbits of the action of $H$ on $V$ are $V_s$ with $s\in S$. $C_e=\{\langle s \rangle e,s\in S\}$, so $C_e$ intersects every $V_s$.

\bigskip\noindent (2) Let $u=\langle s \rangle e\in C$. $Stab_Hu=\{\delta_g|\delta_g(\langle s \rangle e)=\langle s \rangle e\}=\{\delta_{s^n},n\in\ZZ\}=\langle \delta_s \rangle$. 

\bigskip\noindent (3) Let $t \in S$ not necessarily distinct from $s$. Since the elements of $V_{t}$ form a partition of $G$, $\forall n\in\ZZ,\exists! v\in V_{t}:s^n\in v$. Hence $\forall n\in\ZZ,\exists! v\in V_{t}:(uv,s^n)\in E$. Therefore, $\langle \delta_s\rangle$ acts regularly on the set of edges adjacent to $u$ with their other end-point in $V_{t}$.

\bigskip\noindent\textbf{Conversely}, let $\Gamma=(V,E,\psi)$ be a graph with $H$ and $C$ verifying the above property. For all $u\in C$, let $\sigma_u$ be a generator of the cyclic group $Stab_Hu$ (if $C$ is infinite, this requires the axiom of choice). Let $S$ be the multiset $S=\{\sigma_u, u\in C\}$. We will prove that $\Gamma\cong\Psi(H,S)$.

\bigskip\noindent Let $\Psi(H,S)=(V',E',\psi')$, and let $f:V'\rightarrow V$ be the map $\langle \sigma_u\rangle h\mapsto h^{-1}(u)$. For each pair $(u,v)\in C^2$, we choose an edge $(uv,x_0)$ from the multi-edge $uv$ (again, axiom of choice). Let $f^\#:E'\rightarrow E$ be the map $(\langle \sigma_u\rangle h\langle \sigma_v\rangle h',h'')\mapsto h''^{\#-1}(uv,x_0)$. $(f,f^\#)$ is a homomorphism of graphs, as indeed, if $h''\in \langle \sigma_u\rangle h\cap \langle \sigma_v\rangle h'$, then $h''^{-1}(u)=h^{-1}(u)$ and $h''^{-1}(v)=h'^{-1}(v)$. In order to prove that $(f,f^\#)$ is an isomorphism, we now build its inverse $(g,g^\#)$.

\bigskip\noindent Let $u\in C$ and let $O_u$ be its orbit. Supposing that $C$ contains another element $v$ in $O_u$, the set of edges adjacent to $u$ with their other end-point in $O_u$ would contain a link $a$ between $u$ and $v$ as well as a loop $b$ on $v$. However, $Stab_Hu$ acts regularly on this set of edges. Hence, there should be a shift that sends loop $b$ to link $a$ which is impossible. Therefore $C$ contains only one element of $O_u$.

\bigskip\noindent Let $w$ be a vertex in $V$. There exists a unique $u\in C$ such that $u$ is in the orbit of $w$.  Since $u$ is in the orbit of $w$, there exists $h\in H$ such that $h(w)=u$. $h$ is not necessarily unique, however if $\exists (h,h')\in H^2$ such that $h(w)=h'(w)=u$, then $hh'^{-1}\in Stab_Hu$ and therefore $\langle \sigma_u\rangle h=\langle \sigma_u\rangle h'$. Hence $\langle \sigma_u\rangle h$ is independent of the choice of $h$ and entirely determined by $w$. Therefore we can define the map $g:V\rightarrow V'$, $w\mapsto \langle \sigma_u\rangle h$.

\bigskip\noindent Let $l\in E$, with end-points $w'$ and $w'$. There exists a unique couple $(u,v)\in C^2$ such that $u\in O_{w}$ and $v\in O_{w'}$. Also, there exists a couple $(h,h')\in H^2$ with $h(w)=u$ and $h'(w')=v$. Now, $h(w)=u$ and $h(w')\in O_{w'}$. Hence $h^\#(l)$ is adjacent to $u$ and has its other end-point in $O_{w'}$, like $(uv,x_0)$. So, since $Stab_Hu$ acts regularly on the set of edges adjacent to $v$ with their other end-point in $O_{w'}$, there exists a unique $h''\in Stab_Hu$ such that $h''^\#(h^\#(l))=(uv,x_0)$. $h''h(w')=v$, so $h''hh'^{-1}(v)=v$. Hence $h''h\in\langle \sigma_u\rangle h \cap \langle \sigma_v\rangle h'$. Therefore, there is an edge $(\langle \sigma_u\rangle h\langle \sigma_v\rangle h',h''h)$ in $E'$.

\bigskip\noindent We can now define the map $g^\#:E\rightarrow E', l\mapsto (\langle \sigma_u\rangle h \langle \sigma_v\rangle h',h''h)$. If $w$ and $w'$ are the end-points of $l$, then $g(w)$ and $g(w')$ are the end-points of $g^\#(l)$. So $(g,g^\#)$ is a homomorphism of graphs. As $(g,g^\#)$ is the inverse of $(f,f^\#)$, $\Gamma\cong\Psi(H,S)$.
\qed

\begin{corollary}\label{Stab}
Let $\Gamma=\Phi(G,S)$ be a $G$-graph and let $H$ be the set of shifts of $\Gamma$. For all $u\in V(\Gamma)$, for all $s\in S$, $Stab_Hu$ is cyclic and acts regularly on the set of edges adjacent to $u$ with their other extremity in $V_u$.
\end{corollary}

\noindent\textbf{Proof:} We have already proven that $H$ and $C_e$ verify the conditions of the characterisation of $G$-graphs. Since the colour clique $C_g$ is the image of $C_e$ by the automorphism $\delta_g$, $H$ and $C_g$ also verify the same conditions, for all $g\in G$.

\begin{corollary}\label{characterisation2} Characterisation of $G$-graphs \\
A multigraph without loops $\Gamma=(V,E,\psi)$ is a $G$-graph if and only if there is a subgroup $H$ of $Aut(\Gamma)$ and a clique $C$ such that:
\begin{itemize}
\item each orbit of the action of $H$ on $V$ is a stable
\item $C$ intersects every orbit of the action of $H$ on $V$
\item $\forall u\in C, Stab_Hu$ is cyclic
\item $\forall u\in C$, for each orbit $O$ of $H$ such that $u\notin O$, $Stab_Hu$ acts regularly on the set of edges adjacent to $u$ with their other end-point in $O$.
\end{itemize}
\end{corollary}

\noindent\textbf{Proof:}
Suppose that $\Gamma, H,$ and $C$ verify the above. $\forall v\in C, \forall u$ in the orbit of $v$, add a set of loops in bijection with $Stab_Hu$. The result verifies the previous characterisation.

\section{Complete Bipartite $G$-Graphs}

\noindent A {\em connected component} of a graph is a maximal set of vertices linked by finite walks.

\bigskip\noindent A multigraph $\Gamma=(V,E,\psi)$ is said to be {\em bipartite} if $V$ is the union of two stables. Equivalently, a multigraph is bipartite if it has a proper 2-colouring. A {\em complete} graph is a graph in which any two vertices are adjacent. A {\em complete bipartite} graph is a bipartite graph in which any two vertices, if taken one in each stable, are adjacent.

\begin{proposition}\label{components}
Let $\Gamma=\Phi(G,S)$. Each connected component of $\Gamma$ is isomorphic to $\Phi(\langle S\rangle,S)$, and the set of connected components of $\Gamma$ is in one-to-one correspondence with the set of right cosets of $\langle S\rangle$.
\end{proposition}

\noindent\textbf{Proof:}
The first point may be found in \cite{orbitgraphs}, Cor 1.

%It has been shown in \cite{bre2,bre3} that $\Phi(G,S)$ is connected if and only if $S$ generates $G$. Hence $\Phi(\langle S\rangle,S)$ is connected.

%\bigskip\noindent Let $\Gamma_0$ be the connected component of $\Gamma$ containing the colour clique of $e$, and let $\langle s\rangle x$ be a vertex of $\Gamma_0$. There is a finite walk that links $\langle s\rangle e$ to $\langle s\rangle x$. Also, if two vertices $\langle t\rangle y$ and $\langle u\rangle z$ are adjacent, then $yz^{-1}\in\langle t,u\rangle\subset \langle S\rangle$. Hence $x\in \langle S\rangle$, and therefore, the canonical injection $\Phi(\langle S\rangle,S)\hookrightarrow \Phi(G,S)$ induces an isomorphism between $\Phi(\langle S\rangle,S)$ and $\Gamma_0$.

%\bigskip\noindent Now let $\langle s\rangle x$ be any vertex of $\Gamma$. The shift $\delta_x$ is an automorphism of $\Gamma$ that sends $\Gamma_0$ to the connected component of $\langle s\rangle x$. Hence each connected component of $\Gamma$ is isomorphic to $\Phi(\langle S\rangle ,S)$.

\bigskip\noindent The edges of the connected component $\Gamma_0$ containing the colour clique of $e$ are labeled with the elements of $\langle S\rangle$. Since any connected component of $\Gamma$ is the image of $\Gamma_0$ by a shift, it follows that its edges are labeled with the elements of a right coset of $\langle S\rangle$ in $G$. Hence, each connected component of $\Gamma$ corresponds to a right coset of $G$.

\bigskip\noindent An {\em abelian $G$-graph} is a $G$-graph $\Phi(G,S)$ where $G$ is an abelian group.

\begin{proposition}\label{components2}
If $\Gamma_0$ is a connected $G$-graph (resp. abelian $G$-graph), and $\Gamma$ is a graph with all its connected components isomorphic to $\Gamma_0$, then $\Gamma$ is a $G$-graph (resp. abelian $G$-graph).
\end{proposition}

\noindent\textbf{Proof:} Let $\Gamma_0=\Phi(G,S)$ and let $X$ be the set of connected components of $\Gamma$. Let us assume that $X$ is in bijection with an abelian group $(Y,\cdot ,e_Y)$. (The assertion stating that every set is in bijection with an abelian group is a consequence of the axiom of choice.) It stems from proposition~\ref{components} that $\Gamma=\Phi(G\times Y,\{(s,e_Y),s\in S\})$.

\begin{proposition}\label{bipartite}
Let $(s,t)\in G^2,s\neq t$. $\Phi(\langle s,t \rangle ,\{s,t\})$ is complete bipartite if and only if $\forall x\in\langle s,t \rangle,\exists m,n\in\ZZ:x=s^mt^n$.
\end{proposition}

\noindent\textbf{Proof:} $\Phi(\langle s,t\rangle,\{s,t\})$ is complete bipartite \\
$\eq \forall x,y\in\langle s,t \rangle,\langle s \rangle x\cap\langle t \rangle y\neq\emptyset \\
\eq \forall x,y\in\langle s,t \rangle,\exists m,n\in\ZZ:s^mx=t^ny \\
\eq \forall x,y\in\langle s,t \rangle,\exists m,n\in\ZZ:xy^{-1}=s^mt^n \\
\eq \forall x\in\langle s,t \rangle,\exists m,n\in\ZZ:x=s^mt^n$

\begin{proposition}\label{abelian}
Let $\Gamma=\Phi(G,S)$ be an abelian $G$-graph. For any $(s,t)\in S^2, s\neq t$, all the connected components of the induced sub-graph $\Gamma[V_s\cup V_t]$ are isomorphic and are complete bipartite.
\end{proposition}

\noindent\textbf{Proof:} The sub-graph induced by the levels of $s$ and $t$ is isomorphic to $\Phi(G,\{s,t\})$. It stems from the previous two propositions that each of its components is isomorphic to $\Phi(\langle s,t \rangle,\{s,t\})$ which is complete bipartite.\qed

\bigskip\noindent The following result is a generalisation of Prop 2.7 in \cite{gilibert}.

\bigskip\noindent Let $K_{m,n}^l$ be the complete bipartite multi-graph without loops with $m$ vertices on one level, $n$ on the other, and with all its multi-edges of multiplicity $l$, $(l,m,n)\in\mathbb(N^*)^3$.

\begin{proposition}\label{Klmn}
$K_{m,n}^l$ is an abelian $G$-graph.
\end{proposition}

\noindent\textbf{Proof:} For any prime number $p$, and for $n\in\mathbb{N}$, let $v_p(n)$ denote the valuation of $p$ in the prime decomposition of $n$. Let us split the prime factors of $l$ into two groups, depending on whether their valuation is greater in $m$ or in $n$: let $I$ be the set of prime factors $p$ of $l$ such that $v_p(m)\geq v_p(n)$, and let $J$ be the set of prime factors $q$ of $l$ which verify $v_q(m)< v_q(n)$.

\bigskip\noindent Let $l_1=\prod\limits_{p\in I} p^{v_p(l)}$, $l_2=\prod\limits_{q\in J} q^{v_q(l)}$, $d_1=\prod\limits_{q\in J} q^{v_q(m)}$ and $d_2=\prod\limits_{p\in I} p^{v_p(n)}$.\\
Thus, $l=l_1l_2$, $d_1|m\wedge n $ and $d_2|m\wedge n$.\\
Let $G=(\mathbb{Z}/ml_1\mathbb{Z}) \times (\mathbb{Z}/nl_2\mathbb{Z})$, let $s=(1,\frac n{d_1})$ and $t=(\frac m {d_2},1)$, and let $S=\{s,t\}$. We will prove that $K_{m,n}^l\cong\Phi(G,S)$.

\bigskip\noindent Note that in this proof $G$ is an abelian group with an operation noted additively. Hence $\langle S\rangle$ is the set of linear combinations of $s$ and $t$ with coefficients in $\ZZ$. In order to show that $\Phi(G,S)$ is connected, we need to show that $S$ is a generating set.

\bigskip\noindent $\frac m{d_2}s-t=(0,\frac{mn}{d_1d_2}-1)\in\langle S\rangle$. $\frac{mn}{d_1d_2}-1$ and $nl_2$ are coprime. Hence, $(0,1)$ is a multiple of $\frac m{d_2}s-t$, and $(0,1)\in\langle S\rangle$. So $(1,0)=s-\frac n{d_1}\times (0,1)\in\langle S\rangle$, and $S$ is a generating set. Since $G$ is abelian, it now stems directly from proposition~\ref{abelian} that $\Phi(G,S)$ is complete bipartite.

\bigskip\noindent Let $k\in\ZZ$. \\
$ks=(0,0)\eq ml_1|k$ and $nl_2|k\frac n{d_1}$ \\
$\eq \exists k':k=k'ml_1$ and $nl_2|k'\frac{mnl_1}{d_1}$\\
$\eq \exists k':k=k'ml_1$ and $l_2|k'$ \\
$\eq ml_1l_2|k$\\
$\eq ml|k$ \\
Hence $o(s)=ml$, and similarly, $o(t)=nl$.

\bigskip\noindent The elements of $V_s$ form a partition of $G$ into cosets of cardinal $o(s)$. Hence $|G|=o(s)|V_s|$, and $|V_s|=\frac{|G|}{o(s)}=\frac{mnl}{ml}=n$. Similarly $|V_t|=m$. 
It stems from proposition~\ref{structure} that each multi-edge of $\Phi(G,S)$ is in one-to-one correspondence with a right coset of $\langle s\rangle\cap\langle t\rangle$ via its labeling. Hence, the multiplicity of each multi-edge is equal to $|\langle s\rangle\cap\langle t\rangle|=\frac{o(s)o(t)}{|G|}=l$. Therefore, $\Phi(G,S)\cong K_{m,n}^l$.

\section{Incidence Graphs}

The question of finding out when the incidence graph of a simple orbit graph is also a simple orbit graph was asked in \cite{orbitgraphs}. In this section, after recalling some basic properties, we consider two extreme cases of this problem: the bipartite case and the case of complete graphs. These correspond respectively to the cases where $S$ is a pair and where $G$ is a singlet.

\subsection{Definition and basic properties}

\noindent Let $\Gamma=(V,E,\psi)$ be a multigraph. The {\em incidence graph} of $\Gamma$ is the simple graph $I\Gamma$ such that $V(I\Gamma)=V\cup E$ and for each edge $a$ of $\Gamma$, there are two edges in $E(I\Gamma)$ linking $a$ and its end-points.

\bigskip\noindent The following result is a generalisation of Lemma 3 in \cite{balcor}.

\begin{proposition}\label{canonicalinjection}
Given a multigraph $\Gamma$, there is a canonical injection $Aut(\Gamma)\hookrightarrow Aut(I\Gamma),f\mapsto \tilde f$. The image of the canonical injection is the subgroup of $Aut(I\Gamma)$ of automorphisms which stabilise both $V$ and $E$.
\end{proposition}

\noindent\textbf{Proof:} Let $f\in Aut(\Gamma)$. Set $\tilde f(v)=f(v)$ for all $v\in V$, and $\tilde f(a)=f^\#(a)$ for all $a\in E$.  If $v$ and $a$ are adjacent in $I\Gamma$, then $v$ is an end-point of $a$. Hence $\tilde f(v)$ and $\tilde f(a)$ are adjacent. As $I\Gamma$ is simple, $\tilde f^\#$ is entirely determined by $\tilde f$. Therefore $\tilde f$ is a well defined automorphism of $I\Gamma$.

\bigskip\noindent Injectivity: if $\tilde f=\tilde f'$, then $\forall v\in V, f(v)=\tilde f(v)=\tilde f'(v)=f'(v)$ and $\forall a\in E, f^\#(a)=\tilde f(a)=\tilde f'(a)=f'^\#(a)$. Hence $f=f'$. $f\mapsto\tilde f$ is injective.

\bigskip\noindent Image of $f\mapsto\tilde f$: by construction, $\tilde f$ stabilises both $V$ and $E$. Conversely, let $\phi\in Aut(I\Gamma)$ be an automorphism which stabilises both $V$ and $E$. Set $f(v)=\phi(v)$ for all $v\in V$, and $f^\#(a)=\phi(a)$ for all $a\in E$. $f$ is an automorphism of $\Gamma$ and $\phi=\tilde f$.

\begin{proposition}\label{simplebipartite}
Let $\Gamma=\Phi(G,\{s,t\})$ be a simple bipartite $G$-graph. If $o(t)=2$, then $\Gamma$ is an incidence graph.
\end{proposition}

\noindent\textbf{Proof:} Define the graph $\Gamma'=(V_s,V_t,\psi_{\Gamma'})$ with $\psi_{\Gamma'}(\langle t\rangle x)=\{\langle s\rangle x,\langle s\rangle tx\}$. $V(I\Gamma')=V_s\cup V_t=V(\Gamma)$. $V_s$ and $V_t$ are stables in both $I\Gamma'$ and $\Gamma$.

\bigskip\noindent Also, $\langle s\rangle x$ and $\langle t\rangle y$ are adjacent in $I\Gamma' \eq \langle s\rangle x$ is an end-point of $\langle t\rangle y$ in $\Gamma'$ \\
$\eq \langle s\rangle x\in \psi_{\Gamma'}(\langle t\rangle y)=\{\langle s\rangle y,\langle s\rangle ty\}$ \\
$\eq \exists m: s^mx\in\{y,ty\}$ \\
$\eq \langle s\rangle x\cap\langle t\rangle y\neq\emptyset$ \\
$\eq \langle s\rangle x$ and $\langle t\rangle y$ are adjacent in $\Gamma$

\bigskip\noindent Hence, since $\Gamma$ and $I\Gamma'$ are simple graphs with the same vertices and adjacency relations, $\Gamma\cong I\Gamma'$.

\bigskip\noindent Remark: It stems from proposition~\ref{structure} that $\Phi(G,S)$ is a simple graph if and only if $s$ and $t$ are independent.

\subsection{Incidence Graphs of bipartite $G$-graphs}

\noindent The following result is a generalisation of Prop 1 in \cite{orbitgraphs} and Th 3 in \cite{balcor}.

\begin{theorem}\label{incidencebipartite} %{th:incidence $G$-graph}
Let $\Gamma=\Phi(G,S)$ with $S=\{s,t\}$. If $I\Gamma$ is a $G$-graph, then there exists an application $f:\langle s,t\rangle \rightarrow\langle s,t\rangle $ such that: $f$ is involutive, $f(e)=e$ and $\forall x\in\langle s,t\rangle ,\exists m,n\in\ZZ:f(sx)=t^mf(x)$ and $f(tx)=s^nf(x)$. Conversely, if there is a homomorphism $f:\langle s,t\rangle \rightarrow\langle s,t\rangle $ such that: $f$ is involutive, $f(e)=e$ and $\forall x\in\langle s,t\rangle ,\exists m,n\in\ZZ:f(sx)=t^mf(x)$ and $f(tx)=s^nf(x)$, then $I\Gamma$ is a $G$-graph.
\end{theorem}

\noindent Note that $f$ is a homomorphism only in the sufficient condition.

\bigskip\noindent\textbf{Proof:}
It stems from propositions~\ref{components} that the connected components of $\Gamma$ are all isomorphic to $\Phi(\langle S \rangle,S)$. Hence, each connected component of $I\Gamma$ is isomorphic to $I\Phi(\langle S \rangle,S)$. It now stems from proposition ~\ref{components2} that $I\Gamma$ is a $G$-graph if and only if $I\Phi(\langle S \rangle,S)$ is a $G$-graph. Hence we need only consider the case where $G=\langle s,t\rangle $, i.e. where $\Gamma$ is connected.

\bigskip\noindent Suppose that $\Gamma$ is connected and $I\Gamma$ is a $G$-graph: $I\Gamma=\Phi(G',S')$. $I\Gamma$ is bipartite, so $|S'|=2$. To each edge of $\Gamma$ corresponds a vertex of degree 2 in $I\Gamma$. Hence $S'$ contains an element of order 2.

\bigskip\noindent Let $S'=\{s',t'\}$, where $o(t')=2$. $I\Gamma=\Phi(G',S')$, so $V(I\Gamma)=V(\Gamma)\cup E(\Gamma)=V_{s'}\cup V_{t'}$. Suppose that $\langle t'\rangle x'\in E(\Gamma)$. We will show that $V_{s'}=V(\Gamma)$ and $V_{t'}=E(\Gamma)$. Since $\Gamma$ is connected, each vertex (resp. edge) of $\Gamma$ is connected in $I\Gamma$ to $\langle t'\rangle x'$ by a finite walk of even (resp. odd) length. Therefore $V(\Gamma)\subset V_{s'}$ and $E(\Gamma)\subset V_{t'}$. Hence, $V_{s'}=V(\Gamma)$ and $V_{t'}=E(\Gamma)$.

\bigskip\noindent Let $H'$ be the set of shifts of $I\Gamma$. Let $u$ be the edge of $\Gamma$ labeled $e$, which joins $\langle s\rangle e$ and $\langle t\rangle e$. It stems from corollary~\ref{Stab} that $Stab_{H'}u$ is cyclic of order 2. Let $\tilde\tau$ be its generator. Since $\tilde\tau$ is a shift of $I\Gamma$, it stabilises both levels of $I\Gamma$. Hence, it stems from proposition~\ref{canonicalinjection} that $\tilde\tau$ is the image of an automorphism $\tau$ of $\Gamma$ by the canonical injection $Aut(\Gamma)\hookrightarrow Aut(I\Gamma)$.

\bigskip\noindent Since $|SI=2$, it stems from proposition~\ref{structure}(4) that for any $x\in G$, there is exactly one edge labeled $x$ in $\Gamma$. Thus we can define an application $f:G\rightarrow G$ such that $\tau$ sends the edge labeled $x$ to the edge labeled $f(x)$. Since $\tilde\tau$ is of order 2, $f$ is involutive. Also, $Stab_{H'}u$ acts regularly on the set edges adjacent to $u$ (see corollary~\ref{Stab}). Hence $\tau(\langle s\rangle e)=\langle t\rangle e$ and $\tau(\langle t\rangle e)=\langle s\rangle e$. Therefore $f(e)=e$.

\bigskip\noindent Since $\Gamma$ is connected, for all $x\in G$, there is a finite walk that links $\langle s\rangle e$ to $\langle s\rangle x$. And since $\Gamma$ is bipartite and $\langle s\rangle e$ and $\langle s\rangle x$ are on the same level, this walk is of even length. Its image by $\tau$ is a finite walk of even length linking $\langle t\rangle e$ and $\tau(\langle s\rangle x)$. Hence $\tau(\langle s\rangle x)\in V_t$. Similarly, $\tau(\langle t\rangle x)\in V_s$. Therefore, $\tau$ switches $V_s$ and $V_t$.

\bigskip\noindent The edges labeled $x$ and $sx$ are adjacent and meet in $\langle s\rangle x$. Hence their images, labeled $f(x)$ and $f(sx)$ respectively, are adjacent and meet in $V_t$. Therefore, there is a vertex in $V_t$ which contains both $f(x)$ and $f(sx)$. Hence, $\exists m\in\ZZ:f(sx)=t^mf(x)$. Similarly, $\exists n\in\ZZ:f(tx)=s^nf(x)$.

\bigskip\noindent {\bf Conversely}, suppose that there is a homomorphism $f:\langle s,t\rangle \rightarrow\langle s,t\rangle $ such that: $f$ is involutive, $f(e)=e$ and $\forall x\in\langle s,t\rangle ,\exists m,n\in\ZZ:f(sx)=t^mf(x)$ and $f(tx)=s^nf(x)$.

\bigskip\noindent Let $\tau$ be the automorphism of $\Gamma$ which sends the edge labeled $x$ to the edge labeled $f(x)$, for all $x\in G$. $\tau$ is well defined since if the edges labeled $x$ and $y$ are adjacent and meet in $V_s$, then $\exists m\in\ZZ:y=s^mx$, hence $\exists n\in\ZZ:f(y)=t^nf(x)$. So $f(x)$ and $f(y)$ are also adjacent. And similarly if $x$ and $y$ meet in $V_t$, then $f(x)$ and $f(y)$ meet in $V_s$. Therefore $\tau$ preserves adjacency and switches $V_s$ and $V_t$.

\bigskip\noindent Let $H=\langle\tau,\delta_s\rangle$ and let $\widetilde H$ be the image of $H$ by the canonical injection $Aut(\Gamma)\hookrightarrow Aut(I\Gamma)$. Let $u$ be the edge of $\Gamma$ labeled $e$, let $v=\langle s\rangle e$, and let $C=\{u,v\}$. We will show that $H$ and $C$ satisfy the conditions of the characterisation of $G$-graphs.

\bigskip\noindent Let $x\in G$. Since $G=\langle s,t\rangle $, $x$ can be written as $x=\prod\limits_{k=1}^n t^{i_k}s^{j_k}$. Let $\alpha\in\NN$ be such that $s^\alpha=f(t^{i_n})$. $f(f(xs^{-j_n})s^{-\alpha})=\prod\limits_{k=1}^{n-1} t^{i_k}s^{j_k}$. Hence, by induction, there exists integers $\alpha_1,...,\alpha_{2n}$ such that $\prod\limits_{k=1}^{2n} f\circ m_s^{\alpha_k}(x)=e$, where $m_s$ is the right multiplication by $s$. $\prod\limits_{k=1}^{2n} \tau\circ \delta_s^{\alpha_k}$ sends the edge labeled $x$ to $v$. Hence, $H$ acts transitively on $E(\Gamma)$. $\prod\limits_{k=1}^{2n} \tau\circ \delta_s^{\alpha_k}(\langle s\rangle x)=v$ and $\tau\circ\prod\limits_{k=1}^{2n} \tau\circ \delta_s^{\alpha_k}(\langle t\rangle x)=v$. Hence $H$ acts transitively on $V(\Gamma)$. Therefore, the action of $\widetilde H$ on $V(I\Gamma)$ has two orbits: $V(\Gamma)$ and $E(\Gamma)$.

\bigskip\noindent These orbits are stables, and $C$ intersects both of them.

\bigskip\noindent Let $g=\prod\limits_{k=1}^n \tau\circ \delta_s^{\alpha_k}\in Stab_Hu$. Since $g(u)=u$, $\prod\limits_{k=1}^n f\circ m_s^{\alpha_k}(e)=e$. Let $x\in G$.
$\prod\limits_{k=1}^n f\circ m_s^{\alpha_k}(x)=\prod\limits_{k=1}^{n-1} f\circ m_s^{\alpha_k}(f(x)f(s^{\alpha_n}))=...=f^n(x)\prod\limits_{k=1}^{n} f\circ m_s^{\alpha_k}(e)=f^n(x)$
Hence $g=Id$ or $\tau$. Therefore, $Stab_{\widetilde H}u=\{Id,\tilde\tau\}$, which is cyclic and acts regularly on the set of edges adjacent to $u$.

\bigskip\noindent Now, let $g=\prod\limits_{k=1}^n \tau\circ \delta_s^{\alpha_k}\in Stab_Hv$. $g(\langle t\rangle e)$ is adjacent to $v$, so there exists $n$ such that $g(\langle t\rangle e)=\langle t\rangle s^n$. Thus $\delta_s^{-n}g\in Stab_Hu\cap Stab_Hv=\{Id\}$. Hence $g\in \langle\delta_s\rangle$, and $Stab_Hv=\langle\delta_s\rangle$, which is cyclic and acts regularly on the set of edges adjacent to $v$.

\bigskip\noindent Thus $I\Gamma$ verifies the characterisation of $G$-graphs.

\subsection{Incidence Graphs of Complete Graphs}

\noindent Let $K_n$ denote the complete simple graph with $n$ vertices. In this section, we try to find out for which values of $n$ the incidence graph of $K_n$, noted $IK_n$, is a $G$-graph.

\bigskip\noindent Let $n\in\NN, n\geq 2$. Let $\rho$ be the element of the symmetric group $S_n$ defined as follows: $\rho:n-1\mapsto n, n\mapsto n-1$ and $k\mapsto n-1-k$, for $1\leq k\leq n-2$. Note that $\rho^2=Id$. Let $\sigma$ be the following cycle of order $n-1$: $\sigma=(1,2,...,n-1)$.

\begin{theorem}\label{incidencecomplete}
$IK_n$ is a $G$-graph if and only if there exists a permutation $\tau\in S_n$ of order $2$ such that $\tau(n)=n-1$ and $\forall k\in\{1,...,n-2\}, \tau\sigma^k\tau=\sigma^{\tau(k)}\tau\sigma^{\tau\rho\tau(k)}$. In that case, $IK_n\cong\Phi(\langle\sigma,\tau\rangle,\{\sigma,\tau\})$.
\end{theorem}

\noindent\textbf{Proof:}
Suppose that there exists $G$ and $S$ such that $IK_n\cong\Phi(G,S)$. $IK_n$ is connected, bipartite and contains vertices of degree $n-1$ and vertices of degree $2$. Hence $S$ is a generating set with two elements of order $n-1$ and $2$ respectively.

\bigskip\noindent Proposition~\ref{structure} states that $G\simeq \{\delta_g,g\in G\}$ and that each shift $\delta_g$ stabilises both levels of $\Phi(G,S)$. It stems from proposition~\ref{canonicalinjection} that $\delta_g$ is the image of an automorphism of $K_n$ via the canonical injection $Aut K_n\hookrightarrow Aut IK_n$. Finally, by numbering the vertices of $K_n$, we obtain an isomorphism between $Aut K_n$ and $S_n$. Hence $G$ is isomorphic to a subgroup of $S_n$. (Note that while any finite group is isomorphic to a subgroup of {\em a} symmetric group, it is non-trivial that $G$ is isomorphic to a subgroup of the symmetric group {\em of order $n$}.) Therefore, there exists $(s,t)\in S_n^2$ such that $o(s)=n-1, o(t)=2$ and $IK_n\cong\Phi(\langle s, t\rangle,\{s,t\})$.

\bigskip\noindent Recall that every permutation can be decomposed into cycles and that conjugacy classes of $S_n$ are characterised by their cyclic structure. Let us consider the decomposition of $s$ in cycles: suppose that $s$ contains a cycle of length $m>1$. Since $s^m$ has at least $m$ fixed points, $\delta_{s^m}$ stabilises at least $m$ vertices. Since $IK_n$ is simple, $\delta_{s^m}$ also fixes the edges between those vertices. However, if a shift sends an edge to an edge with the same label, it can only be $\delta_{Id}$. So $s^m=Id$ and $m=n-1$. Therefore, $s$ is a cycle of order $n-1$. $s$ and $\sigma$ have the same cycle structure, hence they are conjugate. 

\bigskip\noindent It is shown in \cite{breCayley} that a group isomorphism from $(G_1,S_1)$ to $(G_2,S_2)$ induces a graph isomorphism $\Phi(G_1,S_1)\cong\Phi(G_2,S_2)$. In particular, if $\sigma$ and $\tau$ are the images of $s$ and $t$ respectively by the same conjugation, then $\Phi(\langle\sigma,\tau\rangle,\{\sigma,\tau\})\cong\Phi(\langle s,t\rangle,\{s,t\})$. We will use a first conjugation to replace $s$ with $\sigma$ and a second one to set $\tau(n)=n-1$.

\bigskip\noindent Since $s$ and $\sigma$ are conjugate, there exists $g\in S_n$ such that $gsg^{-1}=\sigma$. Let $a=gtg^{-1}(n)$. If $\sigma$ and $gtg^{-1}$ both fix $n$, then every element of $\langle\sigma,gtg^{-1}\rangle$ fixes the vertex of $IK_n$ corresponding to the vertex of $K_n$ numbered $n$. This contradicts the fact that $\langle\sigma,gtg^{-1}\rangle$ acts transitively on the level of $\sigma$ in $\Phi(\langle\sigma,gtg^{-1}\rangle,\{\sigma,gtg^{-1}\})$. Hence $a\neq n$.

\bigskip\noindent Let $\tau=\sigma^{-a}gtg^{-1}\sigma^a$. Since $\sigma^{-a}gsg^{-1}\sigma^a=\sigma$, $\sigma$ and $\tau$ are the images of $s$ and $t$ via the conjugation by $\sigma^{-a}g$. Hence $\Phi(\langle\sigma,\tau\rangle,\{\sigma,\tau\})\cong IK_n$. Also, $\tau(n)=\sigma^{-a}gtg^{-1}\sigma^a(n)=n-1$.

\bigskip\noindent $\tau$ is a permutation of order 2, so $\tau\neq Id$. Hence $\tau$ does not stabilise any edge, so $\tau$ has at most one fixed point. Therefore, the cyclic decomposition of $\tau$ comprises $\lfloor\frac n 2\rfloor$ transpositions, and $\tau$ has a fixed point if and only if $n$ is odd.

\bigskip\noindent Let $k\in\{1,...,n-2\}$. Since $\tau$ switches $n$ and $n-1$, $\tau(k)\notin\{n-1,n\}$. Let $\phi(k)=n-1-\tau(k)$. $\phi(k)\in\{1,...,n-2\}$ and $\sigma^{\phi(k)}\tau\sigma^k(n-1)=\sigma^{\phi(k)}\tau(k)=n-1$. Therefore, $\exists\phi:\{1,...,n-2\}\rightarrow\{1,...,n-2\}$ such that $\forall k\in\{1,...,n-2\},\sigma^{\phi(k)}\tau\sigma^k\in Stab_G(n-1)$.

\bigskip\noindent However, $Stab_G(n-1)=\langle\tau\sigma\tau\rangle$. Hence $\{\sigma^{\phi(k)}\tau\sigma^k, 1\leq k\leq n-2\}=\{\tau\sigma^k\tau, 1\leq k\leq n-2\}$ ; i.e. $\exists \psi,\chi:\{1,...,n-2\}\rightarrow\{1,...,n-2\}$ such that $\forall k\in\{1,...,n-2\},\tau\sigma^k\tau=\sigma^{\psi(k)}\tau\sigma^{\chi(k)}$.

\bigskip\noindent $\tau\sigma^k\tau(n)=\sigma^{\psi(k)}\tau\sigma^{\chi(k)}(n)$ and $\tau\sigma^k\tau(n-1)=\sigma^{\psi(k)}\tau\sigma^{\chi(k)}(n-1)$. Hence $\tau(k)=\psi(k)$ and $n-1=\psi(k)+\tau(\chi(k))$. So $\chi(k)=\tau(n-1-\psi(k))=\tau(n-1-\tau(k))=\tau\rho\tau(k)$. Thus $\forall k\in\{1,...,n-2\}, \tau\sigma^k\tau=\sigma^{\tau(k)}\tau\sigma^{\tau\rho\tau(k)}$

\bigskip\noindent {\bf Conversely:} suppose that there exists $\tau\in S_n$ of order 2 such that $\tau(n)=n-1$ and $\forall k\in\{1,...,n-2\}, \tau\sigma^k\tau=\sigma^{\tau(k)}\tau\sigma^{\tau\rho\tau(k)}$. Since $\tau(n)\neq n$, $\tau\notin\langle\sigma\rangle$. Hence $\sigma$ and $\tau$ are independent. It then stems from proposition~\ref{simplebipartite} that $\Phi(\langle\sigma,\tau\rangle,\{\sigma,\tau\})$ is the incidence graph of a graph $\Gamma$.

\bigskip\noindent Suppose that $\Phi(\langle\sigma,\tau\rangle,\{\sigma,\tau\})$ contains a $4$-cycle $(\langle\tau\rangle w,a_1,\langle\sigma\rangle x,a_2,\langle\tau\rangle y,a_3,\langle\sigma\rangle z,a_4,\langle\tau\rangle w)$. Then $\langle\tau\rangle w\cap\langle\sigma\rangle x$, $\langle\sigma\rangle x\cap\langle\tau\rangle y$, $\langle\tau\rangle y\cap\langle\sigma\rangle z$ and $\langle\sigma\rangle z\cap\langle\tau\rangle w$ are non-empty. Hence, there exists $(i_1,...,i_4,j_1,...,j_4)$ such that: $\tau^{i_1}w=\sigma^{j_1}x$, $\sigma^{i_2}x=\tau^{j_2}y$, $\tau^{i_3}y=\sigma^{j_3}z$ and $\sigma^{i_4}z=\tau^{j_4}w$.\\
Hence $\tau^{i_1}w=\sigma^{j_1}x=\sigma^{j_1-i_2}\tau^{j_2}y=\sigma^{j_1-i_2}\tau^{j_2-i_3}\sigma^{j_3}z=\sigma^{j_1-i_2}\tau^{j_2-i_3}\sigma^{j_3-i_4}\tau^{j_4}w$\\
And $\sigma^{j_1-i_2}\tau^{j_2-i_3}\sigma^{j_3-i_4}\tau^{j_4-i_1}=Id$. 

\bigskip\noindent If $\sigma^{j_1-i_2}=Id$, then $\tau^{i_1}w=\sigma^{j_1}x=\sigma^{i_2}x=\tau^{j_2}y$, and $(\tau)w=(\tau)y$, which contradicts the assumption that $(\langle\tau\rangle w,a_1,\langle\sigma\rangle x,a_2,\langle\tau\rangle y,a_3,\langle\sigma\rangle z,a_4,\langle\tau\rangle w)$ is a cycle. So $\sigma^{j_1-i_2}\neq Id$. Similarly, $\sigma^{j_3-i_4}\neq Id$ and $\tau^{j_2-i_3}=\tau^{j_4-i_1}=\tau$.
However, the fixed points of $\sigma$ and $\tau\sigma\tau$ are $n$ and $n-1$ respectively. Hence, if $\sigma^j\tau\sigma^k\tau=Id$, then $j$ and $k$ are both multiples of $n-1$. 
Therefore, there are no $4$-cycles in $\Phi(\langle\sigma,\tau\rangle,\{\sigma,\tau\})$. Hence $\Gamma$ is a simple graph.

\bigskip\noindent It stems from the condition $\forall k\in\{1,...,n-2\}, \tau\sigma^k\tau=\sigma^{\tau(k)}\tau\sigma^{\tau\rho\tau(k)}$ that every element of $\langle\sigma,\tau\rangle$ can be expressed using $\tau$ once at most: every element of $\langle\sigma,\tau\rangle$ is of the form $\sigma^k$ or $\sigma^j\tau\sigma^k$, with $0\leq j,k\leq n-2$. Hence $|\langle\sigma,\tau\rangle|\leq  n(n-1)$. In order to prove the equality, we must show that these permutations are all distinct. As established in the previous paragraph, if $\sigma^j\tau\sigma^k\tau=Id$, then $j$ and $k$ are both multiples of $n-1$.  Hence if $(j,k)\neq(l,m)$, then $\sigma^j\tau\sigma^k\neq\sigma^l\tau\sigma^m$. Also, since $\sigma^j\tau\sigma^k(n)\neq n$, $\sigma^j\tau\sigma^k\neq\sigma^i$. Thus $|\langle\sigma,\tau\rangle|=n-1+(n-1)^2=n(n-1)$. Hence $\Phi(\langle\sigma,\tau\rangle,\{\sigma,\tau\})$ has $n(n-1)$ edges. Therefore, $\Gamma$ has $\frac {n(n-1)} 2$ edges.

\bigskip\noindent The number of vertices of $\Gamma$ is $|V_\sigma|=\frac{|\langle\sigma,\tau\rangle|}{o(\sigma)}=n$. So $\Gamma$ is a simple graph with $n$ vertices and $\frac {n(n-1)} 2$ edges. Hence $\Gamma\cong K_n$ and $IK_n\cong\Phi(\langle\sigma,\tau\rangle,\{\sigma,\tau\})$
\qed

\bigskip\noindent Using this result, one may easily check that $IK_n\cong\Phi(\langle\sigma,\tau\rangle,\{\sigma,\tau\})$ for the following values of $n$, and $\tau$:\\
$n=2, \tau=(1,2)$ \\
$n=3, \tau=(1)(2,3)$ \\
$n=4, \tau=(1,2)(3,4)$ \\
$n=5, \tau=(1)(2,3)(4,5)$ \\
$n=7, \tau=(2)(1,5)(3,4)(6,7)$ \\
$n=8, \tau=(1,3)(2,6)(4,5)(7,8)$ \\
$n=9, \tau=(4)(1,2)(3,6)(5,7)(8,9)$ \\
$n=11, \tau=(1)(2,4)(3,6)(5,9)(7,8)(10,11)$ \\
$n=13, \tau=(1)(2,10)(3,4)(5,8)(6,11)(7,9)(12,13)$ \\
$n=16, \tau=(1,12)(3,4)(11,14),(2,9)(6,8)(7,13)(5,10)(15,16)$ \\
$n=17, \tau=(14)(2,8)(1,13)(3,12)(4,15)(5,6)(7,10)(9,11)(16,17)$ \\
$n=19, \tau=(1)(9,17)(3,15)(2,7)(4,11)(14,16)(5,8)(6,10)(12,13)(18,19)$

\bigskip\noindent However, for $n\in\{6,10,12,14,15,18\}$, $IK_n$ is not a $G$-graph. The following propositions explain why.

\begin{corollary}\label{2mod4}
If $n>2$ and $n\equiv 2$ mod $4$, then $IK_n$ is not a $G$-graph.
\end{corollary}

\noindent\textbf{Proof:} $\sigma$ is a cycle of order $n-1$ and $\tau$ has $\lfloor\frac n 2\rfloor$ transpositions. Hence, if $n\equiv 2$ mod $4$, the signatures of $\sigma$ and $\tau$ are $+1$ and $-1$ respectively, thus having $\tau\sigma^k\tau=\sigma^{\tau(k)}\tau\sigma^{\tau\rho\tau(k)}$ is impossible.

\begin{proposition}\label{orbits}
If $IK_n$ is a $G$-graph, then all the orbits of the action of $\langle\rho,\tau\rangle$ on $\{1,...,n\}$ have six elements except for either one or two orbits with two elements
and, if $n$ is odd, one orbit with either one or three elements.
\end{proposition}

\noindent\textbf{Proof:}
If $IK_n=\Phi(\langle\sigma,\tau\rangle,\{\sigma,\tau\})$, then for $k\in\{1,...,n-2\},\tau\sigma^k\tau=\sigma^{\tau(k)}\tau\sigma^{\tau\rho\tau(k)}$. After multiplying by $\tau$ on the left and $\sigma^{-\tau\rho\tau(k)}$ on the right, we obtain $\sigma^k\tau\sigma^{-\tau\rho\tau(k)}=\tau\sigma^{\tau(k)}\tau$. Now, by replacing $k$ with $\tau(k)$, we obtain $\sigma^{\tau(k)}\tau\sigma^{-\tau\rho(k)}=\tau\sigma^k\tau$. Hence $\sigma^{-\tau\rho(k)}=\sigma^{\tau\rho\tau(k)}$, so $\tau\rho\tau(k)=\rho\tau\rho(k), \forall k\in\{1,...,n-2\}$. And finally, since this also holds for $k\in\{n-1,n\}$, $\tau\rho\tau=\rho\tau\rho$. Thus $\langle\rho,\tau\rangle=\{Id,\rho,\tau,\rho\tau,\tau\rho,\rho\tau\rho\}$. Since $|\langle\rho,\tau\rangle|=6$, the orbits of the action of $\langle\rho,\tau\rangle$ on $\{1,...,n\}$ may have $1,2,3$ or $6$ elements.

\bigskip\noindent $|O_k|=2$ if and only if $\tau(k)=\rho(k)$. $\{n-1,n\}$ is an orbit with two elements. Let $(i,k)\in\{1,...,n-2\}^2$ be such that $i+k\neq n-1$. By applying $\tau\sigma^k\tau=\sigma^{\tau(k)}\tau\sigma^{\rho\tau\rho(k)}$ to $\tau(i)$, we obtain: \\
$\tau(k+i)=\tau(k)+\tau(\rho\tau\rho(k)+\tau(i))$ (mod $n-1$). \\
By symmetry, $\tau(k)+\tau(\rho\tau\rho(k)+\tau(i))=\tau(i)+\tau(\rho\tau\rho(i)+\tau(k))$. Hence, if $\tau(i)=\rho(i)$ and $\tau(k)=\rho(k)$, then $\tau(i)=\tau(k)$, and $i=k$. Therefore, if $(i,k)\in\{1,...,n-2\}^2$ is such that $\tau(i)=\rho(i)$ and $\tau(k)=\rho(k)$, then either $i=k$ or $i=\rho(k)$. So other than $\{n-1,n\}$, there is at most one orbit with two elements.

\bigskip\noindent If $n$ is odd, then $\tau$ has a fixed point $m$ and $\rho$ fixes only $\frac{n-1} 2$. Now $\rho\tau\rho(m)=\tau\rho\tau(m)=\tau\rho(m)$. Hence $\tau\rho(m)=\frac{n-1} 2$, and the orbit of $m$ is $O_m=\{m,\rho(m),\frac{n-1} 2\}$. This orbit is a singlet if and only if $m=\frac{n-1} 2$.

\bigskip\noindent Conversely, suppose that $|O_k|=1$ or $3$. $\tau$ stabilises $O_k$ which has an odd cardinal, so since $\tau$ is an involution, $\tau$ has a fixed point in $O_k$. If $n$ is odd, $\tau$ has a unique fixed point $m$, so $O_m$ is the only orbit of cardinal $1$ or $3$. If $n$ is even, $\tau$ has no fixed point, so there is no orbit of cardinal $1$ or $3$.

\begin{corollary}\label{0mod6}
If $6$ divides $n$, then $IK_n$ is not a $G$-graph.
\end{corollary}

\noindent\textbf{Proof:}
If $n$ is even and $IK_n\cong\Phi(\langle\sigma,\tau\rangle,\{\sigma,\tau\})$, then the action of $\langle\rho,\tau\rangle$ on $\{1,...,n\}$ has either one or two orbits with two elements, and all the other orbits have six elements. Since these orbits form a partition of $\{1,...,n\}$, $n\equiv 2$ or $4$ mod $6$.

\begin{corollary}\label{15,21mod24}
If $n\equiv 15$ or $21$ mod $24$, then $IK_n$ is not a $G$-graph.
\end{corollary}

\noindent\textbf{Proof:}
Suppose that $IK_n\cong\Phi(\langle\sigma,\tau\rangle,\{\sigma,\tau\})$. Let $\pi:\{1,...,n-2\}\rightarrow\ZZ/(n-1)\ZZ$, $k\mapsto k-\tau(k)$ mod $n-1$. $\pi$ is injective. Indeed, if $\pi(k)=\pi(l)$ and $k\neq l$, then $\sigma^{\pi(k)}\tau$ has two fixed points, $k$ and $l$. Hence the shift $\sigma^{-\pi(k)}\tau$ fixes an edge, so it is equal to identity, which is impossible.

\bigskip\noindent If $n$ is even, since $\tau$ has no fixed point, then $Im\pi=\ZZ/(n-1)\ZZ-\{0\}$. \\
If $n$ is odd, since $\forall k,\pi(k)\neq\pi(\tau(k))$, then $Im\pi=\ZZ/(n-1)\ZZ-\{\frac{n-1} 2\}$.

\bigskip\noindent Suppose that $n\equiv 3$ mod $6$: \\
The action of $\langle\rho,\tau\rangle$ on $\{1,...,n\}$ has one orbit with two elements, $\{n-1,n\}$, one single point orbit, $\{\frac{n-1} 2\}$, and all the other orbits have six elements.

\bigskip\noindent Let $O_k$ be an orbit with six elements, $O_k=\{k,\rho(k),\tau(k),\rho\tau(k),\tau\rho(k),\rho\tau\rho(k)\}$. $\forall i\in O_k$, $i$ and $\rho(i)$ have the same parity. Thus, depending on the respective parities of $k$, $\tau(k)$, and $\tau\rho(k)$, $\pi$ may take either zero or four odd values on $O_k$. Hence the number of odd values in $Im\pi$ is a multiple of four.

\bigskip\noindent $Im\pi=\ZZ/(n-1)\ZZ-\{\frac{n-1} 2\}$, so the number of odd values in $Im\pi$ is $2\lfloor\frac{n-1} 4\rfloor$. Hence $\lfloor\frac{n-1} 4\rfloor$ is even, and $n\equiv 3$ or $9$ mod $24$.\qed

\bigskip\noindent For $a\in \{1,...,n-1\}, a\wedge (n-1)=1$, let $m_a$ denote the multiplication by $a$ in $\ZZ/(n-1)\ZZ$. We shall identify $\{1,...,n\}$ with $\ZZ/(n-1)\ZZ\cup\{n\}$ and extend $m_a$ to $\{1,...,n\}$ by setting $m_a(n)=n$.

\begin{proposition}\label{conjugates}
If $\tau\in S_n$ satisfies the conditions of Theorem~\ref{incidencecomplete}, i.e. if $\tau$ is of order $2$ and verifies $\tau(n)=n-1$ and $\forall k\in\{1,...,n-2\}, \tau\sigma^k\tau=\sigma^{\tau(k)}\tau\sigma^{\tau\rho\tau(k)}$, then $\forall a\in \{1,...,n-1\}, a\wedge (n-1)=1$, $m_a\tau m_a^{-1}$ also satisfies those same conditions.
\end{proposition}

\noindent\textbf{Proof:}
Conjugation does not change the order, so $m_a\tau m_a^{-1}$ is also of order $2$. Also, $m_a\tau m_a^{-1}(n)=n-1$. Let $k\in\{1,...,n-2\}$, \\
$m_a\tau m_a^{-1}\sigma^k m_a\tau m_a^{-1}=m_a\tau\sigma^{a^{-1}k}\tau m_a^{-1}=m_a\sigma^{\tau (a^{-1}k)}\tau\sigma^{\tau\rho\tau (a^{-1}k)}m_a^{-1}$\\
$=\sigma^{m_a\tau m_a^{-1}(k)}m_a\tau m_a^{-1}\sigma^{m_a\tau\rho\tau m_a^{-1}(k)}=\sigma^{m_a\tau m_a^{-1}(k)}m_a\tau m_a^{-1}\sigma^{m_a\tau m_a^{-1}\rho m_a\tau m_a^{-1}(k)}$\\
Hence $m_a\tau m_a^{-1}$ satisfies the conditions of Theorem~\ref{incidencecomplete}.

\section{Conclusion}
The main contribution of this paper to the study of $G$-graphs is the extension of previous results to the infinite case. By refining the definition of $G$-graphs, we were able to generalise the characterisation and many other results to a wider class of graphs, thus greatly improving our understanding of the underlying structure of $G$-graphs. Regarding the question of identifying when the incidence graph of a $G$-graph is also a $G$-graph, though we have made significant advances both in the bipartite and in the complete cases, an exact characterisation yet remains to be found.

%\begin{ack}
%The authors would like to thank ....
%\end{ack}

%\bibliographystyle{elsart-harv}
%\bibliography{mybibfile}

% Include the ".bib" file (generated by bibtex) right here.

\end{document}